\documentclass{icmart}
\usepackage{epsfig,color}

\contact[Sylvain.Crovisier@math.u-psud.fr]{Laboratoire de Math\'ematiques d'Orsay,
CNRS - UMR 8628, Universit\'e Paris-Sud 11, Orsay 91405, France.}

\newtheorem{theorem}{Theorem}[section]
\newtheorem{corollary}[theorem]{Corollary}
\newtheorem{lemma}[theorem]{Lemma}

\newtheorem{conjecture}[theorem]{Conjecture}
\newtheorem{problem}{Problem}

\theoremstyle{definition}
\newtheorem{definition}[theorem]{Definition}

\def\vdashv{\vdash\!\dashv}
\newcommand{\Diff}{\operatorname{Diff}}
\newcommand{\Tan}{\operatorname{Tang}}
\newcommand{\Per}{\operatorname{Per}}
\newcommand{\Rec}{\operatorname{Rec}}

\newcommand{\cMS}{\mathcal{MS}}
\newcommand{\cG}{\mathcal{G}}
\newcommand{\cX}{\mathcal{X}}
\newcommand{\cI}{\mathcal{I}}
\newcommand{\cR}{\mathcal{R}}
\newcommand{\cU}{\mathcal{U}}
\newcommand{\cV}{\mathcal{V}}
\newcommand{\cW}{\mathcal{W}}
\newcommand{\cO}{\mathcal{O}}
\newcommand{\cP}{\mathcal{P}}

\newcommand{\ZZ}{\mathbb{Z}}
\newcommand{\RR}{\mathbb{R}}
\newcommand{\CC}{\mathbb{C}}
\newcommand{\NN}{\mathbb{N}}

\title{Dynamics of $C^1$-diffeomorphisms: global description and prospects for classification}

\author[Sylvain Crovisier]
{Sylvain Crovisier}

\begin{document}

\begin{abstract}
We are interested in finding a dense part of the space of $C^1$-diffeomorphisms
which decomposes into open subsets corresponding to different dynamical behaviors: we discuss results and questions in this direction.

In particular we present recent results towards a conjecture by J. Palis:
any system can be approximated either by one which is hyperbolic (and
whose dynamics is well understood) or by one which exhibits a homoclinic bifurcation
(a simple local configuration involving one or two periodic orbits).
\end{abstract}

\begin{classification}
Primary 37C20, 37C50, 37D25, 37D30; Secondary 37C29, 37D05, 37G25.
\end{classification}

\begin{keywords} Differentiable dynamical systems, closing lemma, homoclinic tangency, heterodimensional cycle, partial hyperbolicity, generic dynamics.
\end{keywords}

\maketitle

\section{Introduction}

A differentiable transformation - a diffeomorphism or a flow -
on a manifold defines a dynamical systems: our goal is to describe the long time behavior of its orbits.
In some cases, the dynamics, though rich, can be satisfactorily understood:
the hyperbolic systems introduced by Anosov and Smale~\cite{anosov,smale}
break down into finitely many transitive pieces, can be coded by a finite alphabet,
admit physical measures which represent the orbit of Lebesgue-almost every point,
are structurally stable...

The dynamics of a particular system may be quite particular and too complicated.
One will instead consider a large class of systems on a fixed compact connected smooth
manifold $M$ without boundary. For instance:
\begin{itemize}
\item[--] the spaces of $C^r$ diffeomorphisms $\Diff^r(M)$ or vector fields
$\cX^r(M)$, for $r\geq 1$,
\item[--] the subspace $\Diff^r_\omega(M)$ of those
preserving a volume or symplectic form $\omega$,
\item[--] the spaces of $C^r$ Hamiltonians $H\colon M\to \RR$ (when $M$ is symplectic),
and of $C^{r+1}$ Riemannian metrics on $M$ (defining the geodesic flows on
$TM$), etc.
\end{itemize}
This approach (present in~\cite{smale-bourbaki}) allows to study typical dynamics in the class,
but also their stability, i.e. how properties change when the system is replaced by a system nearby.
For finite-dimensional classes of systems (like polynomial automorphisms of $\CC^2$
with fixed degree, or directional flows on flat surfaces with fixed genus)
one can consider sets of parameters with full Lebesgue measure;
for larger classes, one can introduce (non-degenerate)
parametrized families of systems, as in~\cite{palis-yoccoz}.
Working on a Baire space (mainly $\Diff^1(M)$)
we intend here to describe dense subsets of systems that are G$_\delta$
(i.e. Baire-generic) or ultimately even open.

The main difficulty is to \emph{perturb the system while controlling the dynamics}.
Weaker topologies offer more flexibility under perturbations, but less control
on the dynamics. In practice one works in the $C^1$-topology:
for smoother systems, new dynamical properties appear
as Pesin theory~\cite{pesin,pesin-fail}, KAM (see~\cite{yoccoz} and section~\ref{ss.conservative}),
robust homoclinic tangencies on surfaces~\cite{newhouse,moreira},...
but few is known
about perturbations in higher topology
(even about the existence of periodic orbits),
see for instance~\cite{pujals-Cr}.
However, producing $C^1$-open sets, one also describes part of the smoother systems
and presumably gives insights for more regular dynamics.
\medskip

After initial works focused on hyperbolicity, three main problems
have emerged.

\subsection{Density of hyperbolicity}

Smale has explicitly stated~\cite{smale-problem} the following problem
for the class of $C^r$-endomorphisms of the interval (which has been solved affirmatively~\cite{koz-shen-vstrien}) and for the class of one-dimensional complex polynomials with fixed degree
(still unknown).

\begin{problem}[Smale]\label{p.smale}
In which class of systems is hyperbolicity dense?
\end{problem}

In the space of diffeomorphisms $\Diff^r(M)$, $r\geq 1$,
the (open) subset of hyperbolic systems is dense when $\dim(M)=1$,
but this is not the case for any manifold.
Open sets $\cU$ of non-hyperbolic diffeomorphisms have been obtained as follow:
\begin{itemize}
\item[--] When $\dim(M)\geq 3$, Abraham and Smale have built~\cite{abraham-smale} in a non-empty open set $\cU$ a dense family of diffeomorphisms with a {heterodimensional cycle}.
\item[--] When $\dim(M)=2$ and $r\geq 2$, Newhouse has built~\cite{newhouse} in a non-empty open set $\cU$ a dense family of diffeomorphisms with a {homoclinic tangency}.
\end{itemize}
These notions are defined below.
Surprisingly the case $r=1$ and $\dim(M)=2$
is still unknown and has a particular importance for our study.
\smallskip

\begin{figure}[ht]
\vspace{-0.5cm}
\begin{center}
\includegraphics[scale=0.8]{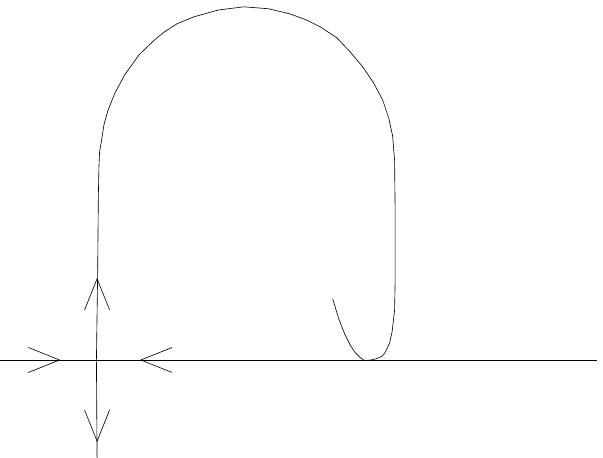} \quad\quad
\includegraphics[scale=0.8]{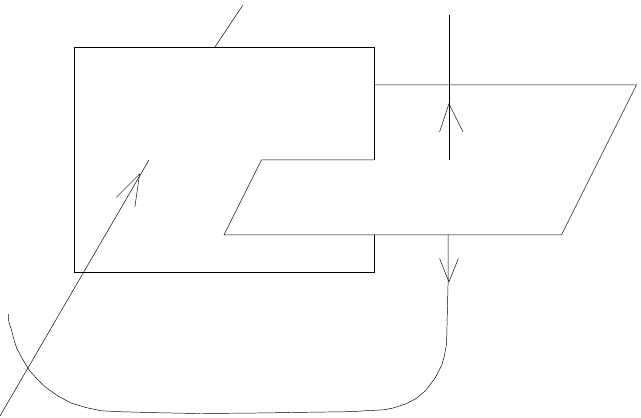}
\begin{picture}(0,0)
\put(-11900,500){$p$}
\put(-4780,2500){$p$}
\put(-1700,2500){$q$}
\end{picture}
\end{center}
\vspace{-0.5cm}
\caption{Homoclinic tangency and heterodimensional cycle.\label{f.bifurcation}}
\end{figure}
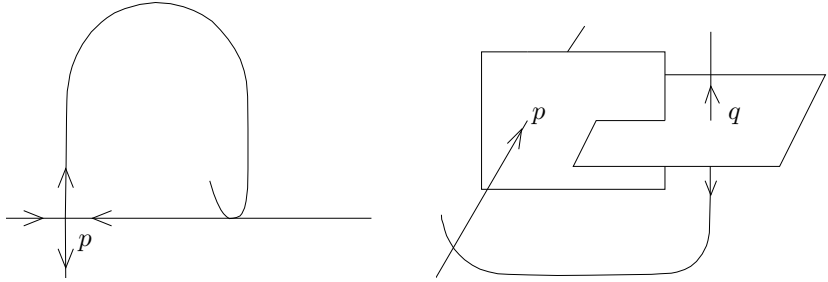

For $f\in \Diff^r(M)$ and any point $p$ in a hyperbolic periodic orbit
(or more generally in a hyperbolic set),
the stable and unstable sets $W^s(p)$ and $W^u(p)$
(i.e. the sets of $z\in M$ such that $d(f^n(p),f^n(z))\to 0$
as $n$ goes to $+\infty$ and $-\infty$ respectively)
are immersed submanifolds with transversal intersection at $p$.
\begin{definition}[Homoclinic bifurcation] \label{d.homoclinic}
A \emph{homoclinic tangency} is a non-trans\-ver\-se intersection
$z\in W^u(p)\cap W^s(p)$ associated to a hyperbolic periodic point $p$.
A \emph{heterodimensional cycle} is a pair of intersections
$z\in W^u(p)\cap W^s(q)$ and $z'\in W^u(q)\cap W^s(p)$
associated to hyperbolic points $p,q$ such that the dimension of $W^s(p)$ is strictly smaller
than the one of $W^s(q)$.
See figure~\ref{f.bifurcation}.
\end{definition}
In both of these configurations the point $z$ is non-wandering and admits a unit tangent vector
whose norm decreases to $0$ under forward \emph{and} backward iterations.

\subsection{Obstructions to hyperbolicity}\label{ss.obstruction}

Palis has conjectured~\cite{palis-conjecture1, palis-conjecture2,palis-conjecture3,palis-conjecture4}
a positive answer to the following problem in the class $\Diff^r(M)$, $r\geq 1$.

\begin{problem}[Palis' conjecture]\label{p.palis}
Approximate any system in a class by
one which is hyperbolic or which exhibits a homoclinic tangency or a heterodimensional cycle.
\end{problem}

The question is to obtain a complete list of simple obstructions to the hyperbolicity.
With Ma\~n\'e's work on stability~\cite{mane-stabilite},
one knows~\cite{hayashi-star,aoki} that a diffeomorphism is non-hyperbolic
if and only if it can be $C^1$-approximated by a diffeomorphism having a non-hyperbolic periodic point.
Two reasons justify that people now look for homoclinic bifurcations
rather than weak periodic orbits.
\smallskip

\noindent
\emph{Cascade of bifurcations, robustness.}
The existence of non-hyperbolic periodic points or homoclinic bifurcations associated
to periodic orbits are one-codimen\-si\-o\-nal configurations and do not occur for open sets of systems.
Replacing the periodic orbits by transitive hyperbolic sets
in definition~\ref{d.homoclinic}, one may obtain open sets
of homoclinic bifurcations and get robust obstructions to hyperbolicity:
this happens for homoclinic tangencies of $C^2$ diffeomorphisms
of surface~\cite{newhouse} and to some extend in higher dimension~\cite{palis-viana};
this also happens in some cases for homoclinic bifurcations in the space of $C^1$ diffeomorphisms when $\dim(M)\geq 3$~\cite{bonatti-diaz-cycle,bonatti-diaz-tangence}.
Indeed for a hyperbolic set, the ``dimension" of its stable set can be larger than the dimension
of its stable leaves.

These homoclinic bifurcations are thus in general not isolated and
as pointed out by Bonatti and D\'{\i}az,
one can strengthen problem~\ref{p.palis}  by requiring the homoclinic tangencies and
heterodimensional cycles to be robust.
\smallskip

\noindent
\emph{Dynamical consequences.}
The unfolding of these bifurcations involve rich dynamics:
homoclinic tangencies generate locally generic sets
of diffeomorphisms displaying infinitely many attracting
or repelling periodic orbits~\cite{newhouse} (which is known as the \emph{Newhouse phenomenon}).
Heterodimensional cycles generate robustly isolated transitive and non-hyperbolic sets~\cite{diaz,
bonatti-diaz-blender}.
\medskip

Pujals and Sambarino have solved Palis conjecture for $C^1$ diffeomorphisms
on surfaces~\cite{ps}.
In higher dimensions some partial results have been obtained, for instance~\cite{wen-palis,pujals-3D1,pujals-3D2,crovisier-modele}.
The following one has been proved in~\cite{crovisier-pujals}.

\begin{theorem}\label{t.crovisier-pujals}
In $\Diff^1(M)$ any diffeomorphism can be approximated
by one which
\begin{itemize}
\item[--] either exhibits a homoclinic tangency or a heterodimensional cycle,
\item[--] or is \emph{essentially hyperbolic}: there exist finitely many hyperbolic attractors
(respectively repellers) whose basin is (open and) dense in $M$.
\end{itemize}
\end{theorem}

In some cases, the dynamics break down into only finitely many pieces, even after perturbation:
these systems, called \emph{tame} are easier to study and may help to test some conjectures.
Using Ma\~n\'e's work, a non-hyperbolic tame dynamics can be perturbed to create
two close periodic points with different stable dimensions. The tameness implies that they
belong to a same piece, hence may be connected in a heterodimensional cycle.
Consequently the Palis conjecture holds in this case for the $C^1$-topology. This includes in particular the conservative
dynamics  (see section~\ref{ss.conservative}).

\subsection{Decomposition of the dynamical space}
\cite{smale-bourbaki,bonatti-survey,crovisier-pujals} propose
to generalize problem~\ref{p.palis} by decomposing (an open and dense subset of)
the considered space of systems into regions which display different dynamical properties.

\begin{problem} Identify new dynamics which allow to split a class of systems.
\end{problem}

\noindent
With Pujals, we suggest~\cite{crovisier-pujals} to focus on two kinds of dynamical properties:
\begin{itemize}
\item[--] \emph{Mechanisms.} We mean simple dynamical configurations which
are non isolated (maybe even robust) and which generate rich dynamical behaviors.
\item[--]\emph{Phenomena.} That is any dynamical property which provides a
global description of the system and holds on a large subset of systems.
\end{itemize}
A mechanism may generate a phenomenon:
for instance the homoclinic tangencies generate the Newhouse phenomenon
for $C^2$ surface diffeomorphisms.
It may also be an obstruction:
one of the first dichotomy was obtained by Newhouse for symplectomorphisms
(hyperbolicity or existence of an elliptic periodic orbit, see theorem~\ref{t.palis-conservatif}).
This mechanism - the existence of an elliptic periodic point - is robust, hence
provides an obstruction to hyperbolicity also for higher topologies.

Another example of dichotomy mechanism/phenomenon is in the following result (which answers a weak version of Palis conjecture). It is proved in~\cite{palis-faible}; the surface and $3$-dimensional cases were obtained before in~\cite{ps} and~\cite{bgw}:

\begin{theorem}\label{t.weak-palis}
The space $\Diff^1(M)$ contains a dense subset $\cMS\cup\cI$
which is the union of two disjoint open sets:
\begin{itemize}
\item[--] $\cMS$ is the set of Morse-Smale diffeomorphisms, i.e. whose
dynamics is hyperbolic and has only finitely many periodic orbits.
Any other orbit accumulates in the future (resp. in the past) towards one of these periodic orbits.
\item[--] $\cI$ is the set of diffeomorphisms $f$ which have a transverse homoclinic orbit:
there exists a hyperbolic periodic point $p$ whose stable and unstable manifolds
have a transverse intersection point different from $p$.\newline
In particular, there exists a compact set $A\subset M$ and an iterate $f^n$
such that $f^n(A)=A$
and the restriction of $f^n$ to $A$ is topologically conjugate to the shift
on $\{0,1\}^\ZZ$. Hence there exists infinitely many periodic orbits.
\end{itemize}
\end{theorem}
The global dynamics for $f\in \cMS$ is very simple, robust under perturbation,
and similar to the time-one map of the gradient flow of a Morse function.
Moreover the topological entropy (which measures the ``complexity" of the system) vanishes.

In the second case, the transverse homoclinic intersection, which is a very simple and robust configuration, implies a very rich behavior, as discovered by Poincar\'e and Birkhoff, and
the topological entropy is non zero. The dynamics however is not described outside a local region
of $M$.
\bigskip

\noindent
{\bf Contents.}
We first discuss generic properties that are consequences of connecting lemmas
for pseudo-orbits. We then present results which led to theorems~\ref{t.crovisier-pujals}
and~\ref{t.weak-palis} above and to other dichotomies inside $\Diff^1(M)$,
see also~\cite{asterisque}.
We  present several questions which emerged during the last years,
and among them some conjectures by Bonatti~\cite{bonatti-survey}.
All these results and questions together constitute a panorama of the main dynamics
which appear in the space of $C^1$-diffeomorphisms.\!\!\!\!\mbox{}
\medskip

Many of the following results were obtained in collaboration with colleagues (and friends),
and in particular with C. Bonatti and E. Pujals.
This subject would have been very different without your viewpoints, thank you!

I am also grateful to F. B\'eguin, R. Potrie, E. Pujals, M. Sambarino, A. Wilkinson and
X. Wang for their comments on the text.

\section{Decomposition of the dynamics}

We say that a dynamical property is $C^1$-generic if it holds
on a dense G$_\delta$ subset of $\Diff^1(M)$.
These properties have been studied in many works, see~\cite{asterisque}.
We describe here properties shared either by \emph{all}
or by \emph{all $C^1$-generic} diffeomorphisms $f$
and that are related to the connecting lemma for pseudo-orbits
stated below.

\subsection{Chain-recurrence}
An open set $U\subset M$ is \emph{attracting} if $f(\overline U)\subset U$.
It decomposes the dynamics into the invariant disjoint compact sets
$A^+=\cap_n f^n(U)$ and $A^-=\cap_n f^n(M\setminus U)$.
The orbits in the complement $M\setminus (A^+\cup A^-)$ are strongly non-recurrent.
By repeating this process one decomposes~\cite{conley} the dynamics into pieces which can be
also obtained with the notion of pseudo-orbits as follows.
\smallskip

For $\varepsilon>0$, a $\varepsilon$-\emph{pseudo-orbit} is a sequence
$(x_n)_{n\in \ZZ}$ such that $d(f(x_n),x_{n+1})<\varepsilon$ for each $n$.
We denote $x\dashv y$ if for each $\varepsilon>0$ there exists a
$\varepsilon$-pseudo-orbit $(x_n)$ and $m\geq 1$ such that $x_0=x$ and $x_m=y$.

\begin{definition}
The \emph{chain-recurrent set} is the (invariant) set $\cR(f)=\{x,\; x\dashv x\}$. 
\end{definition}
\noindent
The chain-recurrent set is compact and contains the set of periodic point $\Per(f)$.
The other classical notions of recurrence - the non-wandering set $\Omega(f)$, the limit set $L(f)$, the recurrent set $\Rec(f)$ - are all contained in $\cR(f)$ and contain $\Per(f)$.\!\!\!\!\mbox{}
\smallskip

On $\cR(f)$ we define the equivalence relation
$x\vdashv y$ whenever there exists a periodic $\varepsilon$-pseudo-orbit
which contains $x,y$ for each $\varepsilon>0$.
\begin{definition}
The chain-recurrent classes of $f$ are the equivalence classes of $\vdashv$.
\end{definition}
\noindent
They are pairwise disjoint invariant compact subsets of the chain-recurrent set.
\smallskip

A chain-recurrence class $\Lambda$ is a \emph{quasi-attractor}
if it admits a basis of attracting neighborhoods.
(It is an \emph{attractor} if $\Lambda=\cap_{n\in \NN} f^n(U)$ for some neighborhood $U$.)
\smallskip

\begin{definition}
Let $K$ be an invariant compact set.
\begin{itemize}
\item[--] $K$ is \emph{chain-transitive} if for any $x,y\in K$
and $\varepsilon>0$ there exists a periodic $\varepsilon$-pseudo-orbit in $K$
which contains $x,y$.
\item[--] $K$ is \emph{transitive} if for any non-empty open sets $U,V$
of $K$, there exists $n\geq 1$ such that $f^n(U)\cap V\neq \emptyset$.
\item[--] $K$ is \emph{topologically mixing} if for any non-empty open sets $U,V$
of $K$, there exists $n_0\geq 1$ such that $f^n(U)\cap V\neq \emptyset$ for each $n\geq n_0$.
\end{itemize}
\end{definition}
\noindent
The chain-recurrence classes are the chain-transitive sets which are maximal for the inclusion.

\begin{definition}
The \emph{homoclinic class} $H(\cO)$ of a hyperbolic periodic orbit $\cO$
is the closure of the transverse intersections between $W^s(\cO)$ and $W^u(\cO)$.
\end{definition}
\noindent
The homoclinic classes satisfy three interesting properties (see~\cite{newhouse-homocline,AC}):
\begin{itemize}
\item[--] $H(\cO)$ contains a dense set of periodic orbits $\cO'$ that are
\emph{homoclinically related} to $\cO$, i.e. such that
$W^u(\cO)$ and $W^s(\cO')$ (resp. $W^u(\cO')$ and $W^s(\cO)$) have a transverse intersection point.
\item[--] $H(\cO)$ is transitive (hence contained in a chain-recurrence class):
there exists a unique $\ell\geq 1$ (called the \emph{period} of $H(\cO)$)
and a subset $A\subset H(\cO)$ such that $f^\ell(A)=A$,
$H(\cO)=A\cup f(A)\dots\cup f^{\ell-1}(A)$ and $f^i(A)\cap A$ has empty interior in
$H(\cO)$ when $0<i<\ell$; moreover $A$ is topologically mixing for $f^\ell$.
\item[--] For any diffeomorphism $g$ that is $C^1$-close to $f$, the orbit $\cO$ has a continuation
$\cO_g$ (given by the implicit function theorem), which gives a notion of
\emph{continuation} $H(\cO_g)$ of a homoclinic class.
\end{itemize}
For general diffeomorphisms, two homoclinic classes may intersect and not coincide.

\subsection{Closing and connecting lemmas in the $C^1$-topology}
In $\Diff^1(M)$ it is possible to perturb one orbit in order to create
periodic points (Pugh's closing lemma~\cite{pugh})
or to connect invariant manifolds of hyperbolic periodic points
(Hayashi's connecting lemma~\cite{hayashi}).
With Bonatti, we have extended~\cite{BC} these technics to pseudo-orbits and obtained:

\begin{theorem}[Connecting lemma for pseudo-orbits]\label{t.connecting}
Let us consider $x,y\in M$ and
assume the following non-resonance condition:
$$\forall n\geq 1,\; \forall v\in TM\setminus \{0\}, \quad Df^n.v\neq v.$$
If $x\dashv y$, there exists
$g$, $C^1$-close to $f$, such that $g^n(x)=y$ for some $n\geq 1$.\\
If $x\in\cR(f)$, there exists
$g$, $C^1$-close to $f$, such that $x$ is periodic for $g$.
\end{theorem}

The perturbation can not be local: one needs to ``close all the jumps"
of a pseudo-orbit. For that purpose we had to build a section of the dynamics:
\begin{lemma}[Topological towers]
There exists $C>0$ (which only depends on $\dim(M)$)
such that for any $f\in \Diff^1(M)$,
any $N\geq 1$ and any (not necessarily invariant) compact set $K$
that does not contain any $i$-periodic point, $1\leq i\leq C.N$,
there exists $U\subset M$ open such that
\begin{itemize}
\item[--] $U$ is disjoint from $f^i(U)$ for $1\leq i\leq N$,
\item[--] $K$ is contained in the union of the $C.N$ first iterates of $U$.
\end{itemize}
\end{lemma}
 
The perturbations in the closing and connecting lemmas may introduce shortcuts in the pseudo-orbits:
for instance theorem~\ref{t.connecting} does not describe the regions which are visited
by the orbit $x,g(x),\dots,g^n(x)$.
Ma\~n\'e~\cite{mane-ergodic} has shown that one can control the distribution
of the periodic orbits in the closing lemma:

\begin{theorem}[Ergodic closing lemma]
There exists a dense G$_\delta$ subset $\cG$ of $\Diff^1(M)$
such that for any $f\in \cG$ and any ergodic probability $\mu$,
there exists a sequence of periodic orbits which converge to $\mu$
for the weak-$*$ topology.
\end{theorem}

The next result~\cite{approximation} gives a topological control on the support of the orbits.

\begin{theorem}[Global connecting lemma]
There exists a dense G$_\delta$ subset $\cG$ of $\Diff^1(M)$
such that any $f\in \cG$ has the following properties:
\begin{itemize}
\item[--] For any points $x_1,\dots,x_k$
satisfying $x_i\dashv x_{i+1}$ for each $1\leq i<k$, and for any $\delta>0$
there exists an orbit of $f$ which intersects each ball $B(x_i,\delta)$.
\item[--] For any chain-transitive set $K$ and any $\delta>0$,
there exists a periodic orbit $\cO$ which is $\delta$-close to $K$
for the Hausdorff topology.
\end{itemize}
\end{theorem}
\medskip

\noindent
{\bf Extension to other classes of systems.}
We stress that the perturbations are supported in a union of small
disjoint balls. This makes difficult the extension of these methods to classes of systems
which do not allow local perturbation.
For the geodesic flow, Contreras has shown~\cite{contreras}
that one can modify the tangent dynamics
above periodic orbits by $C^2$-perturbations, but the following problem is still open:

\begin{problem}
Prove a closing lemma for the geodesic flow (space of $C^2$ metrics).
\end{problem}

\subsection{Decomposition of $C^1$-generic diffeomorphisms}\label{ss.decomposition-generic}
For $C^1$-ge\-ne\-ric diffeomorphisms we have better information on
the chain-recurrence classes.

\paragraph{a- Chain-recurrence classes.}
As a consequence of the connecting and closing lemmas, we obtain
for any $C^1$-generic diffeomorphism $f$:
\begin{itemize}
\item[--] The periodic points are (hyperbolic and) dense in the chain-recurrent set:
$$\overline{\Per(f)}=\overline{L(f)}=\overline{\Rec(f)}= \Omega(f)=\cR(f).$$
\item[--] Any chain-recurrence class
is limit of a sequence of periodic orbits
for the Hausdorff topology.
\item[--] Each chain-recurrence class containing a periodic orbit $\cO$
coincides with the homoclinic class $H(\cO)$.
Any two homoclinic classes are thus disjoint or equal.
\end{itemize}

\paragraph{b- Periodic orbits.}
For $C^1$-generic diffeomorphisms, the periodic orbits inside a homoclinic class $H(\cO)$
have a nice structure (see~\cite{cmp,abcdw,bonatti-diaz-cycle}):
\begin{itemize}
\item[--] Any two periodic orbits in $H(\cO)$ with same stable dimension
are homoclinically related.
\item[--] Any two periodic orbits in $H(\cO)$ with different stable dimension belong to a robust
heterodimensional cycle.
\item[--] The set of stable dimensions of periodic points of $H(\cO)$ is an interval of $\NN$.
\item[--] For any two periodic orbits $\cO_1,\cO_2$ in $H(\cO)$
and any $\theta\in [0,1]$,
there exist periodic orbits in $H(\cO)$ which are arbitrarily close
(for the weak-$*$ topology on finite Borel measures)
to the barycenter $\theta\cdot \cO_1+(1-\theta)\cdot \cO_2$.
\end{itemize}
More about the tangent dynamics above periodic orbits
appear in~\cite{gourmelon-Franks, bochi-bonatti, bcdg,shinohara}.
\paragraph{c- Isolated and tame classes.}
It is equivalent for a chain-recurrence class $\Lambda$ to be isolated
(i.e. $\Lambda$ is open in $\cR(f)$) and to coincide with
the maximal invariant set $\cap_{n\in \ZZ} f^n(U)$ in one of its neighborhoods $U$.
A stronger property is:

\begin{definition}
A chain-recurrence class $\Lambda$ is \emph{tame} if
the maximal invariant set $\Lambda_g:=\cap_{n\in \ZZ} g^n(U)$ in a neighborhood $U$
is a chain-recurrent class for any $g$ $C^1$-close to $f$.
A diffeomorphism is tame if all its chain-recurrence classes are tame.
\end{definition}
If $f$ is $C^1$-generic, an isolated chain-recurrence class is a tame homoclinic class.
Consequently $f$ is tame if and only if it has finitely many chain-recurrence classes.

The chain-recurrence classes of hyperbolic diffeomorphisms
are always isolated (this is part of Smale's spectral theorem).
There are robust examples of isolated chain-recurrence classes
containing periodic points of different stable dimensions
(hence not hyperbolic), see section~\ref{ss.chain-hyperbolic}.
As already noticed at the end of section~\ref{ss.obstruction}, the tame classes
of $C^1$-generic diffeomorphisms are easier to study.
See~\cite{bonatti-diaz-pujals}.\!\!\!\!\!\mbox{}
\medskip

\noindent
\emph{Robust transitivity.}
Let $\Lambda$ be a tame class of  $f$, $C^1$-generic.
One may wonder if $\Lambda_g$ is still transitive or even a homoclinic class
for $g$ $C^1$-close to $f$. A counter example appears in~\cite{bcgp}, but it uses the fact
that it is not a quasi-attractor for $f$ nor for $f^{-1}$.
One can thus ask:
\smallskip

\noindent
\emph{Is any perturbation of $C^1$-generic transitive diffeomorphism still transitive ?}
\smallskip

\noindent
With Abdenur we have answered this question affirmatively assuming
the diffeomorphism is partially hyperbolic with a one-dimensional center bundle.
\medskip

\noindent
\emph{Topological mixing.}
A variation of the connecting lemma gives~\cite{AC} for $C^1$-generic diffeomorphisms
that any isolated homoclinic class $H(\cO)$ decomposes into \emph{disjoint} compact sets
$A\cup f(A)\cup\dots\cup f^{\ell-1}(A)$ where $\ell$ is the period of the class.

In particular one gets the following dichotomy:
\begin{corollary}
There exist two disjoint open sets $\cU_1,\cU_2$ whose union is dense in $\Diff^1(M)$
and which satisfy:
\begin{itemize}
\item[--] $\cU_1$ is the set of diffeomorphisms having a non-empty attracting set $U\neq M$.
\item[--] the diffeomorphisms in a dense G$_\delta$ subset of  $\cU_2$ are topologically mixing.
\end{itemize}
\end{corollary}

\paragraph{d- Non-isolated classes, aperiodic classes.}
Homoclinic classes with robust homoclinic tangencies may create
non-isolated homoclinic classes (accumulated by sinks)
for generic diffeomorphisms:
Newhouse has proved that this occurs on surfaces for the $C^2$-topology,
and in higher dimension for the $C^1$-topology~\cite{newhouse-C1}.

Bonatti and D\'\i az have shown~\cite{bonatti-diaz-aperiodic}
that in some cases these sinks may be turned into non trivial classes:
hence a non-isolated homoclinic class
ejects, after perturbation, new homoclinic classes with similar properties.
This phenomenon is further studied in~\cite{bcdg} and called~\emph{virality}.
Such a $C^1$-generic diffeomorphism present infinite sequences of distinct non-isolated
homoclinic classes, whose limit is a chain-recurrence class disjoint from $\Per(f)$.

\begin{definition}
The chain-recurrence classes which do not contain any periodic point are called \emph{aperiodic classes}.
\end{definition}
\noindent
 Few is known about the dynamics of aperiodic classes:
 the aperiodic classes obtained in~\cite{bonatti-diaz-aperiodic}
 are odometers, but Bonatti and Shinohara are developing a perturbation tool
 which would allow to build non transitive or non uniquely ergodic aperiodic classes.
\smallskip

Some questions remain about non-isolated classes (see also conjectures in~\cite{bonatti-survey}):\!\!\!\!\mbox{}
\smallskip

\noindent
\emph{Is any aperiodic class accumulated by non-isolated (viral) homoclinic classes?}
\smallskip

\noindent
\emph{Is any non-isolated homoclinic class accumulated by aperiodic classes?}
\smallskip

\noindent
One may answer negatively to the second question with
examples of $C^1$-generic diffeomorphisms having no aperiodic classes
and infinitely many homoclinic classes. Indeed, Potrie~\cite{potrie} has built a non-isolated
homoclinic class admitting a neighborhood where the other chain-recurrence
classes are contained in countably many surfaces.
These are homoclinic classes if hyperbolicity is $C^1$-dense on surfaces.

\paragraph{e- Quasi-attractors.}
Theorem~\ref{t.connecting} gives $C^1$-generically:
\begin{itemize}
\item[--] A chain-recurrence class is a quasi-attractor, once
it is \emph{Lyapunov stable}:
there exists a basis of neighborhoods $U$ such that $f(U)\subset U$.
\item[--] There exists a dense G$_\delta$ subset $\cX\subset M$
such that for any $x\in \cX$,
the limit set of the forward orbit $(f^n(x))_{n\geq 0}$
is a quasi-attractor.
\end{itemize}

Attractors may not exist: \cite{bly} gives an example of a $C^1$-generic diffeomorphism
with a quasi-attractor $\Lambda$ which is unique and non-isolated. This quasi-attractor is \emph{essential}:
its basin, i.e. the set of points $x$ such that $f^n(x)$
accumulate on a subset of $\Lambda$ as $n\to +\infty$, is dense in a non-empty open set.
Also the basin of a quasi-attractor may be small:
the aperiodic classes described in~\cite{bonatti-diaz-aperiodic} are quasi-attractors;
each
basin is reduced to the class itself and has empty interior.

\smallskip
One may ask the following for $C^1$-generic diffeomorphisms (see also~\cite{bonatti-survey}):
\smallskip

\noindent
\emph{Is the union of the basins of essential attractors dense in $M$?}
\smallskip

\noindent
\emph{For quasi-attractors is it equivalent to be essential and to be a homoclinic class?}
\smallskip

\noindent
\emph{On attractors, does there exist a \emph{physical measure}?}
(an ergodic probability where the forward orbit of Lebesgue-almost
every point in the basin equidistributes.)

\subsection{Conservative dynamics - ergodicity}\label{ss.conservative}
Conservative systems are chain-transitive. The connecting lemma for pseudo-orbits gives (see~\cite{BC,ABC,AC}):

\begin{theorem}
There exists a dense G$_\delta$ subset $\cG\subset \Diff^1_\omega(M)$
such that any diffeomorphism $f\in \cG$ is topologically mixing.
\end{theorem}
The same statement is false in $\Diff^r_\omega(M)$
when $\omega$ is a volume form and $r$ is large:
by KAM, there may exist a robust one-codimensional invariant torus (see~\cite{yoccoz}).
\smallskip

As already noticed, the $C^1$ Palis conjecture holds in this setting (see~\cite{newhouse-symplectic,palis-faible}):

\begin{theorem}\label{t.palis-conservatif}
In $\Diff^1_\omega(M)$, any diffeomorphism can be approximated
by $f$ which is hyperbolic or which satisfies the following robust property:
\begin{itemize}
\item[] \mbox{}\!\!\!\!\!\!\!\!\!\!\!\!\!\! -- (symplectic case) $f$ has a periodic point with a simple eigenvalue
of modulus $1$.

\item[] \mbox{}\!\!\!\!\!\!\!\!\!\!\!\!\!\! -- (volume case, $\dim(M)\geq 3$)
there exists a robust heterodimensional cycle.
\end{itemize}
\end{theorem}

One can compare to the following~\cite{saghin-xia,horita-tahzibi,dolgopyat-wilkinson,BC}
(see definitions in section~\ref{ss.partial}):

\begin{theorem}
In $\Diff^1_\omega(M)$, any diffeomorphism can be approximated
by one with a completely elliptic periodic point
(eigenvalues are simple, of modulus $1$), or:
\begin{itemize}
\item[] \mbox{}\!\!\!\!\!\!\!\!\!\!\!\!\!\! -- (symplectic case) by one which is partially hyperbolic
and robustly transitive,
\item[] \mbox{}\!\!\!\!\!\!\!\!\!\!\!\!\!\! -- (volume case)
by one which has a (non-trivial) dominated splitting.
\end{itemize}
\end{theorem}
The existence of a completely elliptic periodic point is an obstruction to robust
transitivity~\cite{arbieto-matheus} and
in the symplectic case, the robust transitivity is characterized by partial hyperbolicity.
In the volume-preserving case, Dolgopyat and Wilkinson conjectured~\cite{dolgopyat-wilkinson}:

\begin{conjecture}\label{c.RT}
In the volume preserving case,
the sets of robustly transitive diffeomorphisms and of those having
a dominated splitting have the same closure
in $\Diff^1_\omega(M)$.\!\!
\end{conjecture}
\smallskip

A stronger notion of undecomposability involves
the ergodicity of the volume:
\begin{problem}
Is ergodicity dense (hence Baire-generic) in $\Diff^1_\omega(M)$?
\end{problem}
The $C^1$-generic systems in $\Diff^1_\omega(M)$ with positive metric entropy
are ergodic (this is proved in~\cite{abw} for the symplectic and in~\cite{acw} for the volume preserving cases).
However $C^1$-generic systems with zero metric entropy also occur~\cite{bochi}.
\smallskip

There exists (non-empty) $C^1$-open sets of ergodic diffeomorphisms in $\Diff^r_\omega(M)$
when $r>1$: these diffeomorphisms (which include the hyperbolic systems)
are called \emph{stably ergodic} and were studied intensively
(see~\cite{wilkinson}). Note that it is not known if they exist also in $\Diff^1_\omega(M)$.
Pugh and Shub have conjectured that stable ergodicity is dense in
the space of $C^r$ partially hyperbolic diffeomorphisms.

In parallel to conjecture~\ref{c.RT}, with Avila and Wilkinson we proposed~\cite{acw}:
\begin{conjecture}\label{c.SE}
For $r>1$,
the sets of stably ergodic diffeomorphisms and of those having
a dominated splitting have the same $C^1$-closure
in $\Diff^r_\omega(M)$.
\end{conjecture}
In this direction we obtained~\cite{acw}:
\begin{theorem}
In the space of volume-preserving diffeomorphisms
$\Diff^r_\omega(M)$, $r>1$,
those having a partially hyperbolic splitting $E^s\oplus E^c\oplus E^u$
into non-trivial bundles are contained in the closure of the set of
stably ergodic diffeomorphisms.
\end{theorem}

\section{Notions of weak hyperbolicity}

Let $K$ be an invariant set for $f\in \Diff^1(M)$. We recall the classical notion:
\begin{definition}\label{d.hyperbolicity}
$K$ is (uniformly) \emph{hyperbolic} if
there exists an invariant continuous splitting
$T_K=E^s\oplus E^u$ and $N\geq 1$ such that
$\|Df^N_{|E^s}\|\leq 1/2$ and $\|Df^{-N}_{|E^u}\|\leq 1/2$
(i.e., $E^s$ and $E^u$ are uniformly contracted by $f$ and $f^{-1}$ respectively on $K$).

\noindent
A diffeomorphism is \emph{hyperbolic} if each chain-recurrence class
is hyperbolic.
\end{definition}
It is well-known that hyperbolic sets satisfy several important properties:
they can be continued for diffeomorphisms $C^1$-close,
each of their points has stable and unstable manifolds,
their pseudo-orbits are shadowed by orbits,...
We present now several weaker notions of hyperbolicity which sometimes
keep these properties and will appear in the next sections.

\subsection{Tangent dynamics - partial hyperbolicity}\label{ss.partial}
Pesin theory describes systems where the uniformity in definition~\ref{d.hyperbolicity}
is relaxed: Oseledets theorem associates to any ergodic probability $\mu$
its Lyapunov exponents $\lambda_1\leq\dots\leq \lambda_{\dim(M)}$
which are the possible limits of $\log(\|Df^n(x).u\|)/n$ as $n\to \infty$
for any $u\in T_xM$ and a.e. $x\in M$.
When each $\lambda_i$ is non-zero, $\mu$ is \emph{non-uniformly hyperbolic}.

Here is another way to relax hyperbolicity which allows vanishing exponents.

\begin{definition}
An invariant splitting $T_KM=E\oplus F$ by linear sub-bundles above $K$
is \emph{dominated} if there is $N\geq 1$ such that
$\|Df^N(x).u\|\leq 1/2 \|Df^{N}(x).v\|$ for each $x\in K$
and each unit vectors $u\in E_x$ and $v\in F_x$.
\end{definition}
This definition extends to splittings into a larger number of invariant bundles.
The \emph{finest dominated splitting} is the (unique) one which maximizes this number.

\begin{definition}
A dominated splitting $T_KM=E^s\oplus E^c\oplus E^u$
is \emph{partially hyperbolic} if $E^s$ (resp. $E^u$) is uniformly contracted by $f$ (resp. $f^{-1}$)
and if one of the bundles $E^s,E^u$ is non-trivial.
\end{definition}

Dominated splittings and partial hyperbolicity extend to the closure of $K$
and to invariant sets in a neighborhood of $K$ for diffeomorphisms $C^1$-close.
Moreover any point of a partially hyperbolic set has unique (strong) stable and unstable
manifolds tangent to $E^s$ and $E^u$, that we denote by $W^{ss}(x)$ and $W^{uu}(x)$.

Hirsch, Pugh and Shub have built~\cite{hps} a weak notion of center manifold:

\begin{theorem}\label{t.hps}
If $K$ has a dominated splitting $T_KM=E_1\oplus F\oplus E_2$,
there exists a \emph{locally invariant plaque family tangent to $F$}, i.e. a map
$\cW\colon F\to M$ satisfying:
\begin{itemize}
\item[--] Each induced map $\cW_x\colon F_x\to M$
is an embedding, depends continuously on $x\in K$ for the $C^1$-topology,
$\cW_x(0)=x$, and the image is tangent to $F_x$ at $x$.
\item[--] There exists $\rho>0$ such that $\cW_x(B(0,\rho))$ is sent by $f$ in $\cW_{f(x)}$
for each $x$.
\end{itemize}
\end{theorem}
The image of $\cW_x$ (still denoted $\cW_x$) is a \emph{plaque}.
The plaque family is in general not unique; the union of two different plaques
may not be a sub manifold.

One can sometimes obtain the following stability along the plaques (see~\cite{crovisier-pujals}):
\begin{definition}
The plaque family $\cW$ is \emph{trapped} if $f(\overline{\cW_x})\subset \cW_{f(x)}$ for each $x\in K$.

$F$ is \emph{thin-trapped} if ``inside a plaque family $\cW$",
there exist plaque families which are trapped and whose plaques have arbitrarily small diameters.
\end{definition}

\subsection{Dynamics along one-dimensional center directions}\label{ss.model}
Let $K$ be a chain-transitive set with a dominated splitting $E_1\oplus F\oplus E_2$,
$\dim(F)=1$. Pesin theory does not describe the local dynamics along the direction $F$
when the Lyapunov exponents along $F$ vanish. These dynamics
are studied in~\cite{ps} when $E_1$ or $E_2$ is degenerated.
When they are not, we introduced the next notion~\cite{palis-faible,crovisier-modele}:

\begin{definition}
A \emph{center model} for $K$ is a compact metric space $\widehat K$ and
continuous maps $\widehat f\colon \widehat K\times [0,1]\to \widehat K\times [0,+\infty)$, $\pi\colon \widehat K\times [0,+\infty)\to M$ such that:
\begin{itemize}
\item[--] $\widehat f$ is a local homeomorphism near $\widehat K\times \{0\}$,
\item[--] for each $x$, there is $x'$ satisfying
$\widehat f(\{x\}\times [0,1])\subset \{x'\}\times [0,+\infty)$,
\item[--] $\pi(\widehat K\times \{0\})=K$ and $\pi$ semi conjugates $\widehat f$ and $f$,
\item[--] each $t\mapsto \pi(x,t)$ is a $C^1$-embedding which depends continuously on $x$ for the $C^1$-topology and the image is tangent to $F_{\pi(x,t)}$.
\end{itemize}
\end{definition}
\noindent
From theorem~\ref{t.hps}, $K$ admits
a center model: $\widehat K$ is the unit tangent bundle associated to $E^c$,
hence $\pi$ is two-to-one on $\widehat K\times \{0\}$.

Using pseudo-orbits, the local dynamics of center models can be classified 
into four different types: $\widehat K\times \{0\}$ is a quasi-attractor or not,
for $\widehat f$ or for $\widehat f^{-1}$.
This allows us to prove that one of the following (not exclusive) cases occurs for the local dynamics
along a locally invariant plaque family $\cW$ tangent to $F$:
\begin{description}
\item[Thin-trapped.] \emph{The bundle $F$ is thin-trapped.}

If $E_2$ is uniformly contracted by $f^{-1}$, a weak shadowing lemma implies that the unstable set of $K$ meets
the stable set of a periodic orbit whose chain-recurrence class is non-trivial.
See \cite[Prop. 10.20]{asterisque}, \cite[Prop. 4.5]{crovisier-modele}.

\item[Chain-recurrent.] \emph{In any neighborhood $U$ of $K$,
there exist $x\in K$, a non trivial curve $I$
with $f^n(I)\subset \cW_{f^n(x)}$, for $n\in \ZZ$,
and a chain-transitive set $\Lambda\supset K\cup I$.}

If $U$ is small, any periodic orbit $\cO\subset U$, whose exponent along $F$ is close to $0$
and having a point close to the middle of $I$, belongs to the chain-recurrence class of $K$.
See~\cite[corollary 4.4]{csy}.

\item[Semi chain-unstable.] \emph{There exists a locally invariant
\emph{half} plaque family $\cW^+$.
It is thin-trapped by $f^{-1}$; for any $x\in K$ and $z\in \cW^+_x$ we have $x\dashv z$.}

If $E_1, E_2$ are uniformly contracted by $f,f^{-1}$ respectively,
$f$ is $C^1$-generic and $K$
is not twisted, it is contained in a homoclinic class. See~\cite[Prop. 4.4]{crovisier-modele}.
\end{description}

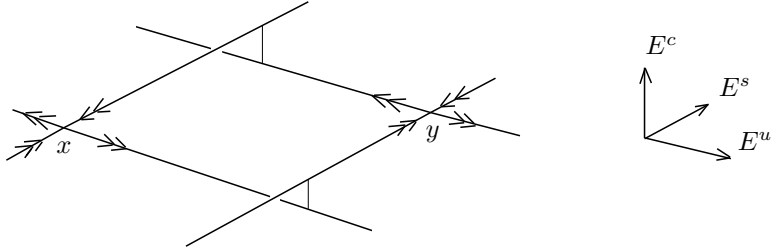
\begin{figure}[ht]
\begin{center}
\includegraphics[scale=0.8]{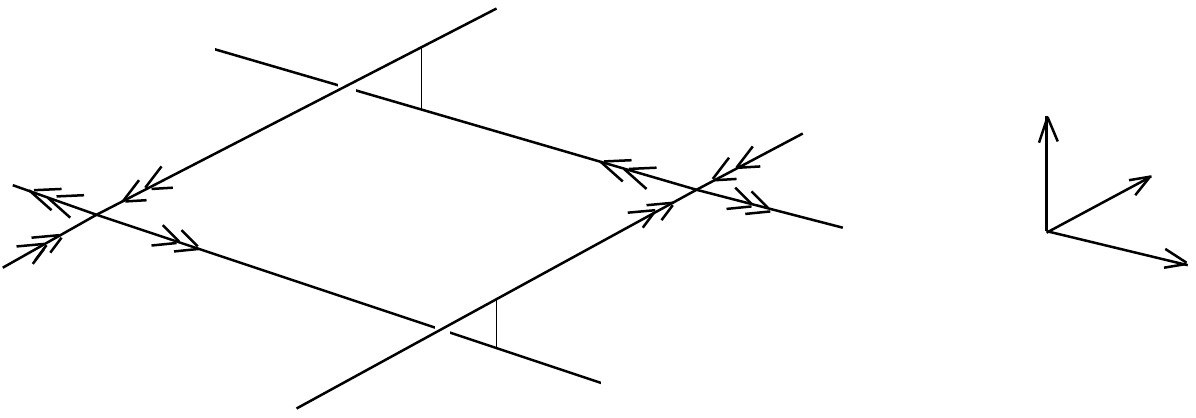}
\begin{picture}(0,0)
\put(-10800,1550){$x$}
\put(-5000,1800){$y$}
\put(-1550,3100){$E^c$}
\put(-100,1600){$E^{u}$}
\put(-400,2450){$E^{s}$}
\end{picture}
\end{center}
\vspace{-0.5cm}
\caption{Two close points in a twisted set. \label{f.twist}}
\end{figure}

The twisted geometry above is very particular.
For the definition,
one extends continuously $F$ in a neighborhood of $K$.
Locally, it is trivial, hence orientable.
\begin{definition}
A partially hyperbolic set $K$ with a one-dimensional center bundle is \emph{twisted}
if for any $x,y\in K$ close, one can connect 
$W^{uu}_{loc}(x)$ to $W^{ss}_{loc}(y)$ and $W^{uu}_{loc}(y)$ to $W^{ss}_{loc}(x)$
by two curves tangent to $F$ having the same orientation.
(Figure~\ref{f.twist}.)
\end{definition}

Using Pugh's and Ma\~n\'e's closing lemma arguments,
when $K$ does not contain periodic point and is twisted,
one can find~\cite[prop. 3.2]{crovisier-modele} a $C^1$-perturbation $g$ having a periodic orbit $\cO$ close to $K$ such that $W^{ss}(\cO)$ and $W^{uu}(\cO)$ intersect.

\subsection{Chain hyperbolicity}\label{ss.chain-hyperbolic}
The following notion is introduced in~\cite{crovisier-pujals}:
\begin{definition}\label{d.chain-hyperbolic}
A homoclinic class $H(\cO)$ is \emph{chain-hyperbolic} if it has a dominated splitting $T_{H(\cO)}M=E^{cs}\oplus E^{cu}$ and plaque families $\cW^{cs}, \cW^{cu}$ such that:
\begin{itemize}
\item[(i)] $\cW^{cs}$ (resp. $\cW^{cu}$) is tangent to $E^{cs}$ (resp. $E^{cu}$)
and trapped by $f$ (resp. $f^{-1}$);
\item[(ii)] there exists $p\in \cO$ such that $\cW^{cs}_p\subset W^s(p)$ and
$\cW^{cu}_p\subset W^u(p)$.
\end{itemize}
$\cW^{cs}$ (resp. $\cW^{cu}$)
are called \emph{center-stable} (resp. \emph{center-unstable}) plaque families.
\end{definition}
\smallskip

\noindent
{\bf Properties.}
At large scale, chain-hyperbolicity looks similar to hyperbolicity and both
notions share several properties:
\begin{itemize}
\item[--] \emph{Robustness.} If $H(\cO)$ is a chain-recurrence class
and $E^{cs},E^{cu}$ are thin-trapped by $f$ and $f^{-1}$, its continuation $H(\cO_g)$ is chain-hyperbolic for $g$ $C^1$-close.
\item[--] \emph{Strong periodic points.} $H(\cO)$ contains a dense set $\cP$ of periodic points $p$
satisfying definition~\ref{d.chain-hyperbolic}(ii) and
whose exponents are bounded away from $0$. They have a hyperbolic continuation for $g$ in a uniform neighborhood of $f$.
\item[--] \emph{Invariant manifolds.} For any $x\in H(\cO)$ and $y\in \cW^{cs}_x$, one has
$y\dashv x$. A transverse intersection of center stable and center unstable plaques still belongs to $H(\cO)$.
\end{itemize}
\smallskip

\noindent
{\bf Examples.}
By deforming a hyperbolic diffeomorphism near a periodic point,
one can build robust examples of non-hyperbolic isolated chain-hyperbolic classes,
see~\cite{bonatti-viana}. One can also build examples from skew product maps.

\smallskip

\noindent
{\it Remark.} Other constructions of robustly transitive sets exist.
The center bundle may be parabolic~\cite{bcgp},
or tangent to a foliation with compact or non-compact leaves~\cite{bonatti-diaz-cycle}.

\section{Dynamics far from homoclinic tangencies...}

If $x$ is a homoclinic tangency for a hyperbolic periodic orbit with
stable dimension $d^s$, there is no dominated splitting $E\oplus F$
with $\dim(E)=d^s$ on the orbit of $x$. We discuss now the converse and
consider the diffeomorphisms $f$ that can not be approximated by
homoclinic tangencies: $f\in \Diff^1(M)\setminus \overline{\Tan}$
where $\Tan$ denotes the collection of diffeomorphisms which exhibit a homoclinic tangency.

Theorems~\ref{t.csy} and~\ref{t.mane} below imply that these dynamics are partially hyperbolic:

\begin{theorem}\label{t.ht}
There exists an open and dense subset $\cU\subset\Diff^1(M)\setminus \overline{\Tan}$ such that
any $f\in \cU$ has at most finitely many sinks and sources
and any of its other chain-recurrence classes $\Lambda$
has a partially hyperbolic splitting:
$$T_\Lambda M=E^s\oplus E^c_1\oplus \dots\oplus E^c_\ell\oplus E^u,$$
where $E^s$ (resp. $E^u$) is non trivial, uniformly contracted (resp. expanded)
and where each $E^c_i$ is one-dimensional.
\end{theorem}
The decomposition of the center into one-dimensional sub-bundles limits the pathological behaviors. In particular, these systems admit symbolic extensions and any continuous
map $\varphi\colon M\to \RR$ has an equilibrium state (see~\cite{LVY,csy}).

\subsection{Existence of weak periodic points inside the class}
Improving a technics of~\cite{ps}, Wen~\cite{wen-tangences}
and Gourmelon~\cite{gourmelon-tangences} have shown:

\begin{theorem}\label{t.gw}
For any diffeomorphism $f\in \Diff^1(M)\setminus \overline{\Tan}$ and any $d^s\geq 1$,
the decomposition into stable and unstable spaces,
above the hyperbolic periodic orbit with stable dimension $d^s$,
is a dominated splitting.
\end{theorem}

With results of section~\ref{ss.decomposition-generic}(a) this implies that for $C^1$-generic
diffeomorphisms far from the homoclinic tangencies,
any chain-recurrence class $\Lambda$ (but maybe finitely many sinks and sources)
has a non-trivial dominated splitting $E\oplus F$.
\smallskip

Assuming that $E$ is not uniformly contracted, we have to build a sequence
of periodic orbit with smaller stable dimension which accumulates on $\Lambda$,
so that $E$ can be further decomposed. Such periodic orbits exist in a neighborhood
of $\Lambda$ (essentially from Ma\~n\'e's ergodic closing lemma), but the difficulty is to
approximate the whole class $\Lambda$ in Hausdorff topology.
There is an easy case: $\Lambda$ contains periodic points with at least one Lyapunov exponent along $E$ which is positive or close to zero.
Indeed in this case $\Lambda$ is a homoclinic class which contains
a dense set of such periodic points and $E$ can be decomposed.

These periodic orbits in $\Lambda$ are obtained under weak hyperbolicity by shadowing:
either from non-uniform hyperbolicity (Liao's ``selecting lemma"~\cite{liao-obstruction2,liao-stability,wen-palis})
or from topological hyperbolicity when the Lyapunov exponents vanish (using center models).
Such weak hyperbolicity fails only when $\Lambda$
has a partially hyperbolic splitting $E^s\oplus E^c\oplus E^u$ with
$\dim(E^c)=1$, $E^s\oplus E^c\subset E$ and such that the Lyapunov exponent along $E^c$ of any invariant probability on $\Lambda$ vanishes ($\Lambda$ is aperiodic).
\smallskip

The stable dimensions of periodic points of $\Lambda$
is an interval in $\NN$ and one gets the following statement proved in~\cite{csy}
(already obtained in~\cite{GWY} for minimal sets):

\begin{theorem}\label{t.csy}
There is a dense G$_\delta$ set $\cG\subset \Diff^1(M)\setminus \overline \Tan$ such that
for $f\in \cG$:\!\!\!\mbox{}
\smallskip

\noindent
Any aperiodic class $\Lambda$ is partially hyperbolic:
$T_\Lambda M=E^s\oplus E^c\oplus E^u$, $\dim(E^c)=1$, and the center Lyapunov exponent
of any invariant probability on $\Lambda$ vanishes.
\smallskip

\noindent
Any homoclinic classes $H(\cO)$ is partially hyperbolic: 
\begin{itemize}
\item[--] $T_{\!H(\cO)\!} M\!\!=\!\!E^s\oplus E^c_1\oplus \dots\oplus E^c_\ell\oplus E^u\!\!$,
\item[--] each $E^c_i$ is one-dimensional and $H(\cO)$ contains (weak) periodic orbits whose Lyapunov exponent along $E^c_i$ is arbitrarily close to $0$;
\item[--] either $E^s\oplus E^c_1$ (resp. $E^c_\ell\oplus E^u$)
is thin-trapped by $f$ (resp. $f^{-1}$) or there is $p\in H(\cO)$ periodic
whose stable (resp. unstable) space is $E^s_p$ (resp. $E^u_p$).
\end{itemize}
\end{theorem}

The last item follows from section~\ref{ss.model}:
there exist periodic points whose stable space contains $E^c_1$, so the semi chain-unstable case
does not occur. If the chain-recurrent one holds, there is a (weak) periodic point $p$ whose
stable space is $E^s_p$.

\subsection{Extreme bundles}
The fact that the uniforms bundles $E^s,E^u$ in theorem~\ref{t.csy} are non trivial
comes from the next result obtained with Pujals and Sambarino~\cite{cps}.
It generalizes Ma\~n\'e's argument for interval endomorphisms~\cite{mane-1D}
and previous works~\cite{ps,ps-homocline,crovisier-pujals}.
This completes the statement of theorem~\ref{t.ht}.

\begin{theorem}\label{t.mane}
For any $f\in \Diff^2(M)$ and any invariant compact set $K$
with a dominated splitting $T_K M=E\oplus F$, $\dim(F)=1$ such that:
\begin{itemize}
\item[--] each periodic point in $K$ has an unstable space containing $F$,
\item[--] there is no periodic closed curve in $K\setminus \Per(f)$ tangent to $F$,
\end{itemize}
 then $F$ is uniformly expanded.
\end{theorem}
An important tool of the proof is the construction of  ``semi-geometrical" Markov rectangles,
that are laminated charts by curves tangent to $F$.

\subsection{Dichotomy Morse-Smale / homoclinic intersections}
For proving theorem~\ref{t.weak-palis},
it is enough to take a $C^1$-generic diffeomorphism $f$
far from homoclinic tangencies whose homoclinic classes are reduced
to isolated periodic orbits (even after perturbation).
One has to consider an aperiodic class.

Any aperiodic class is partially hyperbolic with one-dimensional center
and section~\ref{ss.model} applies.
It can not be twisted since a transverse homoclinic intersection would appear
after perturbation. The three types can be ruled out since they
would give the existence of a non-trivial homoclinic class.
Hence there is no aperiodic class and
theorem~\ref{t.weak-palis} follows.

\subsection{Quasi-attractors}\label{ss.finiteness}
Using the technics of the section~\ref{ss.model}, one gets
more information on quasi-attractors for generic diffeomorphisms
in $\Diff^1(M)\setminus \overline \Tan$:
\begin{itemize}
\item[--] Quasi-attractors are homoclinic classes $H(\cO)$, see~\cite{yang}.
\item[--] Considering the splitting $T_{H(\cO)} M=E^s\oplus E^c_1\oplus \dots\oplus E^c_\ell\oplus E^u$, there exists a periodic orbit in $H(\cO)$ whose unstable dimension is equal to
$\dim(E^u)$, see~\cite[theorem 4]{crovisier-pujals}.
\item[--] Quasi-attractors are essential attractors (proved by\! Bonatti,\! Gan,\! Li, D.\!\! Yang).
\item[--] If moreover $f$ can not be approximated by diffeomorphisms with a heterodimensional cycle,
the number of quasi-attractors is finite, see~\cite{crovisier-pujals}.
\end{itemize}

By studying the geometry of invariant compact sets saturated by the strong unstable leaves,
we proved recently with Sambarino and Potrie the finiteness of the quasi-attractors for $C^1$-generic systems in the class of diffeomorphisms whose chain-recurrences classes are partially hyperbolic with a one-dimensional center bundle. This class offers an ideal setting to study the uniqueness
of physical measures and equilibrium states (see~\cite{rrtu,viana-yang}
for smooth diffeomorphisms and~\cite{qiu} for $C^1$-generic hyperbolic diffeomorphisms).

\subsection{Obstruction to the Newhouse phenomenon}\label{ss.newhouse}
A hyperbolic periodic orbit $\cO$ is \emph{sectionally dissipative}
if its two largest Lyapunov exponents $\lambda_1,\lambda_2$ (counted with multiplicity) satisfy $\lambda_1+\lambda_2<0$.
If such an orbit has a homoclinic tangency, one can obtain a sink by $C^1$-small
perturbation.  Theorem~\ref{t.mane} implies that the converse holds $C^1$-generically:

\begin{corollary}
For any open set $\cV \subset \Diff^1(M)$, the next properties are equivalent:
\begin{itemize}
\item[--] Baire-generic diffeomorphisms in $\cV$ have infinitely many sinks,
\item[--] densely in $\cV$ there exist diffeomorphisms exhibiting homoclinic tangencies associated to sectionally dissipative periodic points.
\end{itemize}
\end{corollary}

One can expect to characterize the absence of Newhouse phenomenon:
\begin{conjecture}
There exist two disjoint open sets $\cU_1,\cU_2$ whose union is dense in $\Diff^1(M)$
and which satisfy the following properties:
\begin{itemize}
\item[--] Baire-generic diffeomorphisms in $\cU_1$ have infinitely many sinks;
\item[--] the diffeomorphisms $f\in\cU_2$ are volume hyperbolic: each chain-recurrence class
$\Lambda$, which is not a sink, has a dominated splitting
$T_\Lambda=E\oplus F$ where $F$ is non-trivial and $|\det(Df^N_{|F})|>1$ on $\Lambda$ for some $N\geq 1$.
\end{itemize}
\end{conjecture}

\section{... and far from heterodimensional cycles}

Theorem~\ref{t.crovisier-pujals}
is now a consequence of the following (from~\cite{crovisier-pujals}):
\begin{theorem}\label{t.cycle}
There exists a dense G$_\delta$ subset $\cG\subset \Diff^1(M)$
such that for any $f\in \cG$, any quasi-attractor which is partially hyperbolic
with a one-dimensional center bundle is either hyperbolic
or contains a (robust) heterodimensional cycle.
\end{theorem}
In this section one considers a quasi-attractor with a partially hyperbolic
splitting $E^s\oplus E^c\oplus E^u$, $\dim(E^c)=1$,
for $f$ $C^1$-generic. By section~\ref{ss.finiteness} it is a homoclinic class $H(\cO)$.
One can assume that all the periodic points
in $H(\cO)$ have stable dimension equal to $\dim(E^s)+1$
and, from theorem~\ref{t.csy}, that $E^s\oplus E^c$ is thin-trapped.

Indeed, by the results of section~\ref{ss.decomposition-generic}(b), the conclusion of the theorem~\ref{t.cycle} holds if
$H(\cO)$ contains two periodic points with different stable dimension
and by section~\ref{ss.finiteness} it always contains periodic points
of stable dimension $\dim(E^s)+1$.

\subsection{Strong homoclinic intersections} Let us assume that there exist
diffeomorphisms $g$ that are $C^1$-close to $f$ and satisfy the following property.

\begin{definition}
$H(\cO_g)$ has a \emph{strong homoclinic intersection} if there exist periodic points
$p,q$ homoclinically related to $\cO_g$ such that $(W^{ss}(p)\setminus \{p\})\cap W^{u}(q)\neq \emptyset$.

An invariant set $K$ with a partially hyperbolic splitting $E^{ss}\oplus F$ has a \emph{strong connection} if it contains a point
$x$ such that $(W^{ss}(x)\setminus \{x\})\cap K\neq \emptyset$.
\end{definition}
By theorem~\ref{t.csy}, if $H(\cO)$ is not hyperbolic it contains weak periodic
points. A strong homoclinic intersection, for a diffeomorphism close,
can be moved on these weak points. This gives a robust heterodimensional cycle
by $C^1$-perturbation.
\smallskip

The non-existence of strong connection allows  to reduce
the dimension of the ambient manifold: the following is a consequence~\cite{whitney} of
Whitney's extension theorem and of a graph transform argument.

\begin{theorem}
Any invariant set $K$ with a partially hyperbolic splitting $E^{ss}\oplus F$
and no strong connection
 is contained in a $C^1$ submanifold $\Sigma$ tangent to $F$
 which is locally invariant: $\Sigma\cap f(\Sigma)$ is a neighborhood of $K$ in $\Sigma$.
\end{theorem}

When $H(\cO)$ is contained in a locally invariant submanifold $\Sigma$ tangent to
$E^c\oplus E^u$, theorem~\ref{t.mane} implies that $E^c$ is uniformly contracted.
Hence we are reduced to the case where a strong connection exists,
i.e. $H(\cO)$ contains $x\neq y$ such that $W^{ss}(x)=W^{ss}(y)$.
\smallskip

In a homoclinic class the periodic points are dense,
hence one can consider $p_x,p_y$ periodic close to $x,y$ respectively
so that $W^{u}_{loc}(p_x)$ and $W^u_{loc}(p_y)$ are close to
the local unstable manifolds of $x,y$.
One can hope that the projections
of the unstable manifolds of $x,y$ by strong stable holonomy are ``topologically transverse".
This implies that there exists $x'\in W^{u}_{loc}(p_x)$, $y'\in W^u_{loc}(p_y)$
such that $W^{ss}(x')=W^{ss}(y')$. Since $H(\cO)$ is a quasi-attractor,
$x',y'$ are still in the class. Other more degenerated cases may occur and have to
be handled by other arguments.
\medskip

We now have to deal with the following problem:
\smallskip

\noindent
{\bf Reduced problem.} {\it Assume that $H(\cO)$ contains periodic points $p_x,p_y$
homoclinically related to $\cO$ and $x\!\neq\! y$ such that $x\!\in\! W^{uu}(p_x)$,
$y\!\in\! W^{uu}(p_y)$, $W^{ss}(x)=W^{ss}(y)$.

Does there exist $g$ near $f$ such that $H(\cO_g)$ has a strong homoclinic
intersection?}

\subsection{Pointwise continuation of chain-hyperbolic classes}
The class $H(\cO)$ is chain-hyperbolic and the results of section~\ref{ss.chain-hyperbolic}
apply.
The points $p_x,p_y$ may be taken in the set of strong periodic points $\cP$
and for any diffeomorphism $g$ in a neighborhood $\cU$ of $f$,
the continuations of the points of $\cP$ are well-defined and dense in $H(\cO_g)$.
This allows to introduce the following notion, similar to the branched holomorphic motion considered by Dujardin and Lyubich in~\cite{dujardin-lyubich}
for holomorphic families of polynomial automorphisms of $\mathbb{C}^2$.

\begin{definition}
For any $g,g'\in \cU$, one says that $x\in H(\cO_g)$ and $x'\in H(\cO_{g'})$
have the \emph{same continuation} if there exists a sequence $(p_n)$ in $\cP$
such that $(p_{n,g})$ converges to $x$ and $(p_{n,g'})$ converges to $x'$.
\end{definition}

A point of $H(\cO_g)$ may have
several
continuations in
$H(\cO_{g'})$.
However in our setting the center-unstable bundle of the chain-hyperbolic structure is uniformly expanded. Consequently for any point  $x$ in the unstable manifold of a point $p_x\in \cP$,
its continuation - denoted by $x_g$ - is unique and depends continuously in $g$.
\medskip

In the setting of the reduced problem above, the continuations $x_g,y_g$ belong
to the same center stable plaques for any $g\in \cU$ and we are led to ask:
\smallskip

{\it Does there exist $g$ near $f$ such that $W^{ss}(x_g)\neq W^{ss}(y_g)$?}
\smallskip

\noindent
Indeed $\cW^{cs}_x\setminus W^{ss}(x)$ has two connected components;
if there exists $g^+,g^-$ such that $y_{g^+}$ and $y_{g^-}$
belong to the continuations of different components of $\cW^{cs}_x\setminus W^{ss}(x)$,
by considering $q_x,q_y\in \cP$ close enough to $x,y$ and an arc
$(g_t)$ in $\cU$ between $g^+$ and $g^-$, one finds a diffeomorphism $g$
such that $W^{ss}(q_{x,g})$ and $W^u(q_{y,g})$ intersect. This gives a strong
homoclinic intersection as required.

\subsection{How to remove a strong connection}
We are still in the (simplest) setting of the reduced problem above and
look for $g\in \cU$ such that the strong connection between $x$ and $y$ is broken.

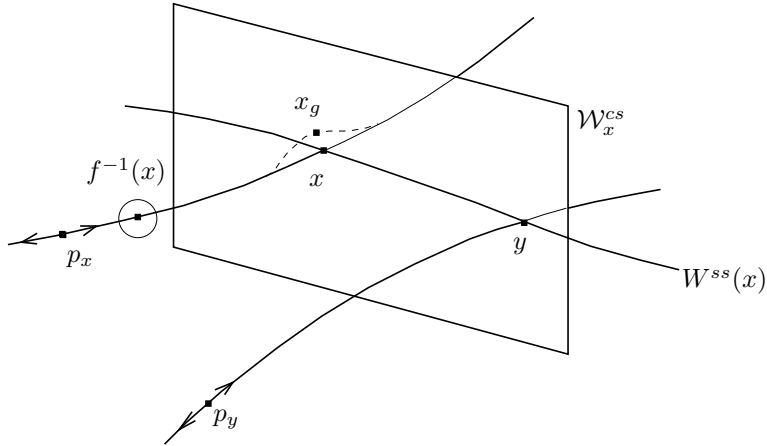
\begin{figure}[ht]
\begin{center}

\includegraphics[width=9cm]{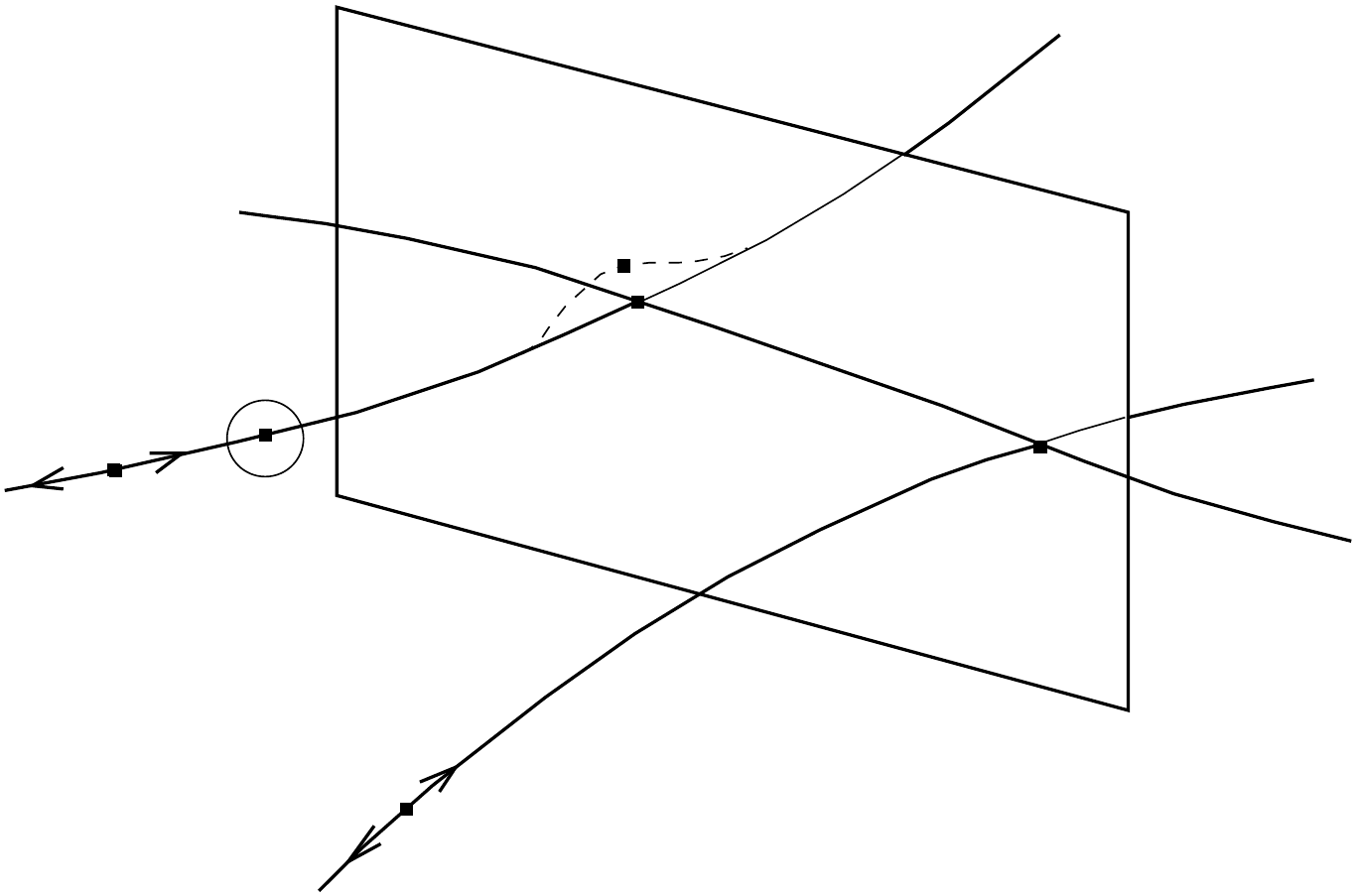}
\begin{picture}(0,0)
\put(-6235,5300){$x_g$}
\put(-6000,4100){$x$}
\put(-2800,3100){$y$}
\put(-1800,5000){$\cW^{cs}_x$}
\put(-150,2500){$W^{ss}(x)$}
\put(-9800,2900){$p_x$}
\put(-7500,400){$p_y$}
\put(-9500,4200){$f^{-1}(x)$}
\end{picture}
\end{center}
\vspace{-0.5cm}
\caption{A broken strong connection.\label{f.connection}}
\end{figure}

The idea is to modify $f$ in a ball $B(f^{-1}(x),r)$
so that $W^u(p_{x,g})$ intersects a given component of
$\cW^{cs}_{x_f}\setminus W^{ss}(x_f)$.
The distance $d(x_g,g(f^{-1}(x)))$ is arbitrarily small with respect to the size $r$
of the support of the perturbation.
If the positive orbit of $y$ does not return ``too fast" in the support of the perturbation,
using the weak hyperbolicity, one shows that
the distances $d(y,y_g)$ and $d(W^{ss}_{loc}(y), W^{ss}_{loc}(y_g))$ are small also
and the connection is broken (figure~\ref{f.connection}).
A different argument is performed when the returns
of the positive orbit of $y$ near $x$ are fast. See~\cite{crovisier-pujals,asterisque}.

\section{Panorama of the dynamics in $\Diff^1(M)$}

We end this text by summing up several questions and conjectures which
allow to structure the space of $C^1$-diffeomorphisms.
Most of them already appear in~\cite{bonatti-survey, asterisque}.

\subsection{Global dynamics}\label{ss.global}
As noticed in~\cite[section 1.3]{bonatti-diaz-cycle},
all the examples of $C^1$-generic non-hyperbolic systems
involve heterodimensional cycles (this becomes false in higher topologies).
This justifies:

\begin{conjecture}[Bonatti-D\'\i az hyperbolicity conjecture]
Any diffeomorphism can be approximated in $\Diff^1(M)$
by one which is hyperbolic or exhibits a robust heterodimensional cycle.
\end{conjecture}

Motivated by the results of sections~\ref{ss.newhouse} and~\ref{ss.finiteness}
we expect a positive answer to following conjecture made by Bonatti in~\cite{bonatti-survey}.
\begin{conjecture}[Bonatti's finiteness conjecture]
In $\Diff^1(M)\setminus \overline{\Tan}$,
there exists an (open and) dense subset of tame diffeomorphisms.
\end{conjecture}

However there exists robust examples of transitive dynamics in $\overline{\Tan}$
and a positive answer to the previous conjecture would not give a dichotomy.

\begin{problem}
Characterize non-tame dynamics: find a robust mechanism which
generates non-tame dynamics and whose union with tame systems is dense in
$\Diff^1(M)$.
\end{problem}

Note that both conjectures imply the $C^1$ Palis conjecture
and that the first one implies a positive answer to Smale problem
for $C^1$-diffeomorphisms on surfaces.
\medskip

\noindent
{\bf A last example: the universal dynamics.}
Let us mention that
there exists a non-empty open set $\cU$ of diffeomorphisms
having a homoclinic class with no dominated splitting, such that the volume is contracted
above one periodic orbit and is expanded above another one.
This implies~\cite{bonatti-diaz-aperiodic} that the dynamics of the $C^1$-generic diffeomorphisms
$f$ in $\cU$ are \emph{universal}: any diffeomorphism $g$ of the unit ball $B(0,1)\subset \RR^{\dim(M)}$
may be approximated by the restriction of some iterates $f^n$ to some balls $B\subset M$.
A similar property holds for $C^\infty$-diffeomorphisms on surfaces~\cite{turaev}.
\medskip

\begin{figure}[ht]
\begin{center}
\vspace{-0.1cm}
\includegraphics[scale=0.7]{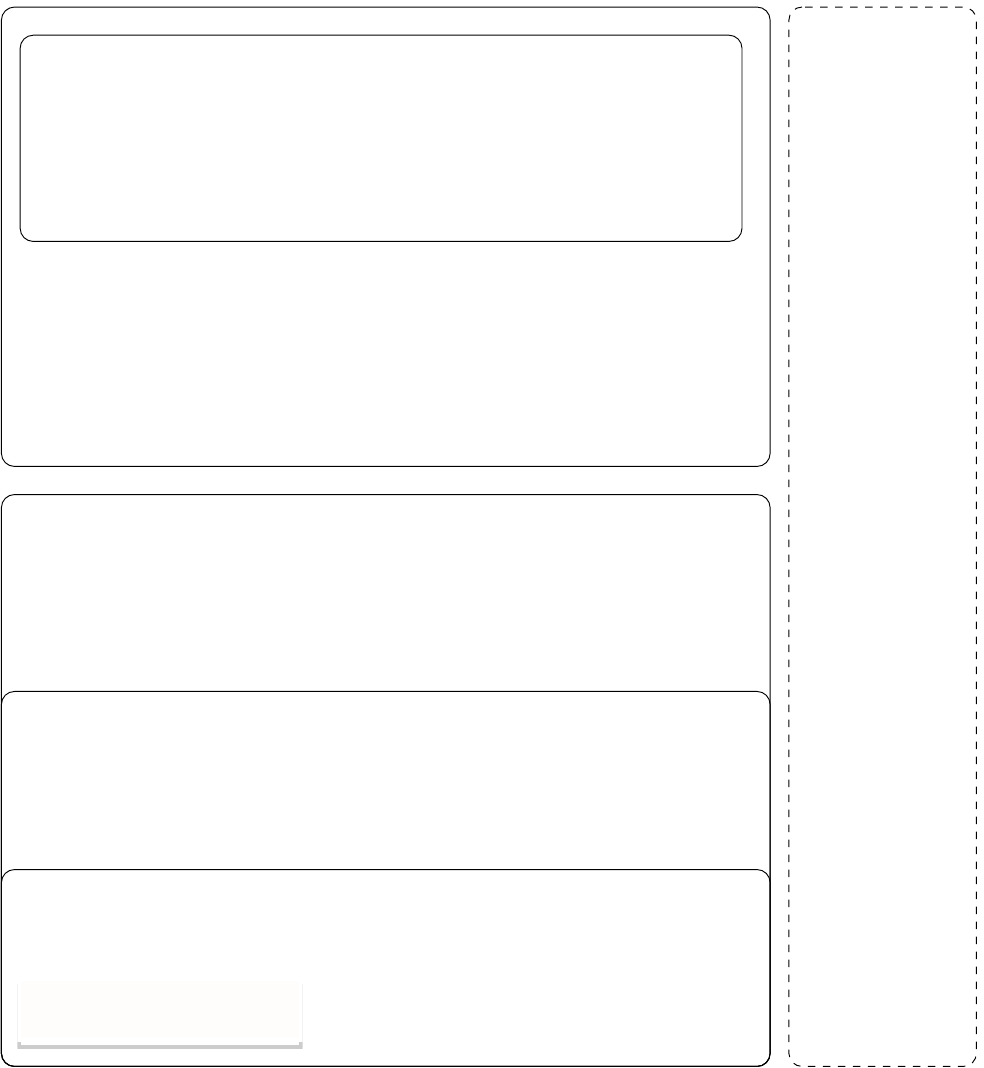}
\begin{picture}(0,0)(6300,-100)
\includegraphics[scale=0.28]{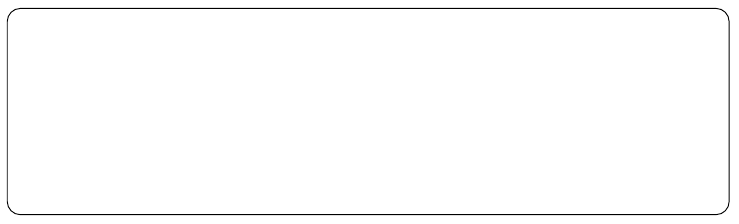}
\end{picture}
\setlength{\unitlength}{1579sp}
\begin{picture}(0,0)
\put(-5330,450){\makebox(0,0){Morse-Smale}}
\put(-5330,1250){\makebox(0,0){hyperbolic}}
\put(-5330,6450){\makebox(0,0){not tame}}
\put(-5330,5950){\makebox(0,0){(critical + heterodimensional)}}
\put(-5330,8000){\makebox(0,0){universal}}
\put(-1450,4000){\rotatebox{90.0}{other?}}
\put(-5330,4450){\makebox(0,0){critical tame}}
\put(-5330,3950){\makebox(0,0){(critical + heterodimensional)}}
\put(-5330,2700){\makebox(0,0){tame heterodimensional}}
\put(-5330,2200){\makebox(0,0){(not critical)}}
\end{picture}
\end{center}
\vspace{-0.5cm}
\caption{Structure of the dynamical space $\Diff^1(M)$.
\label{f.panorama}}
\end{figure}
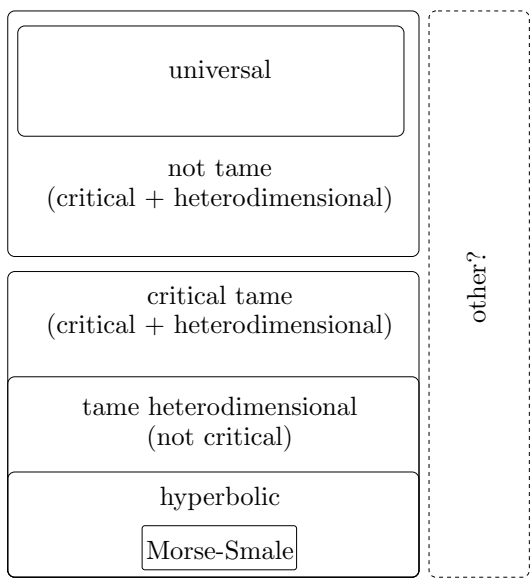

Assuming that the two conjectures above hold,
one can decompose the space $\Diff^1(M)$ into
disjoint regions with increasing complexity, as pictured on figure~\ref{f.panorama}.

\subsection{Local dynamics}\label{ss.local}
The previous conjectures do not control where the homoclinic bifurcations occur.
We state now more precise questions which allow to break
the conjectures into three steps.
\smallskip

Let $f$ be a $C^1$-generic and non-hyperbolic diffeomorphism.

\begin{description}
\item[I. Localization.]
One knows that one of the chain-recurrence classes is non-hyperbolic.
We expect that this is the case for at least one homoclinic class:

\mbox{}\!\!\!\!\!Ia. \emph{Does $f$ exhibit a non-hyperbolic homoclinic class?}

\mbox{}\!\!\!\!\!Ib. \emph{If $f$ is not tame, does it have a non-isolated homoclinic class?}

\item[II. Local dichotomies.] Let $H(\cO)$ be a non-hyperbolic homoclinic class $H(\cO)$.
Strengthening Palis conjecture, one may look for homoclinic bifurcations
inside the class $H(\cO)$ (rather than in a neighborhood).

\mbox{}\!\!\!\!\!IIa. \emph{Does $\cO_g$ belong to a heterodimensional cycle for some $g$ close?}

\mbox{}\!\!\!\!\!IIb. \emph{If $H(\cO)$ is not tame, does $\cO_g$ has a homoclinic tangency for some $g$ close?}

\item[III. Robustness.]
At last, we are aimed to stabilize the homoclinic bifurcation.
This is possible for heterodimensional cycles as shown in~\cite{bonatti-diaz-cycle}.
Let us assume that $\cO_g$ has a homoclinic tangency for some $g$ close to $f$.

\mbox{}\!\!\!\!\!\emph{Does $\cO$ belong to a hyperbolic set with robust homoclinic tangencies?}

(The same question for heterodimensional cycles has been answered in~\cite{bonatti-diaz-cycle}.)
\end{description}

These intermediate questions have been discussed in the case of
Smale's problem for surface diffeomorphisms~\cite{abcd}.
Moreira has shown~\cite{moreira} that robust tangencies do not occur for $C^1$-diffeomorphisms
on surface, solving the step III in this case. A possible approach for the two first steps is to control
the lack of dominated splitting by considering the critical set introduced by Pujals and F. Rodriguez-Hertz~\cite{pujals-hertz}.

\subsection{Tangent dynamics}\label{ss.tangent}
Some of the previous questions may be addressed by a better understanding
of the weak hyperbolicity on each chain-recurrence class.
Considering the known examples, the case of tame diffeomorphisms
or of diffeomorphisms far from the homoclinic tangencies~\cite{bonatti-diaz-pujals,csy},
and the results~\cite{bcdg,cps,whitney}, we formulate the following conjectures.
Recent discussions with X. Wang seem to bring a partial answer towards
the first one.

\begin{conjecture}
Let $H(\cO)$ be a non-hyperbolic homoclinic class for a $C^1$-generic $f$,
and $E^s\oplus E^c_1\oplus\dots\oplus E^c_\ell\oplus E^u$ the
finest dominated splitting
such that $E^s$ (resp. $E^u$) is the maximal uniformly contracted (resp. expanded) sub-bundle.

Then, the minimal stable dimension $k$ of the periodic points satisfies
$$\dim(E^s)\leq k<\dim(E^s\oplus E^c_1).$$
Moreover, when $\dim(E^c_1)\geq 2$ and $\dim(E^s)<k$ two cases are possible:
\begin{itemize}
\item[--] On a periodic orbit, the sum of all the Lyapunov exponents inside $E^c_1$ is positive.
Then $k=\dim(E^s)+1$ and the class is contained in a locally invariant submanifold tangent to $E^c_1\oplus\dots\oplus E^c_\ell\oplus E^u$ (and is not isolated).
\item[--] For any ergodic probability, the sum of the $\dim(E^s\oplus E^c_1)-k+1$ larger Lyapunov exponents inside $E^c_1$ is negative. (The volume along $E^c_1$ is contracted.)
\end{itemize}
\end{conjecture}

The next conjectures implies that for an aperiodic class
with a dominated splitting $E\oplus F$,
either $E$ is uniformly contracted or $F$ is uniformly expanded.
\begin{conjecture}
Let $\Lambda$ be an aperiodic class for a $C^1$-generic $f$,
and $E^s\oplus E^c\oplus E^u$ the
dominated splitting
such that $E^s$ (resp. $E^u$) is the maximal uniformly contracted (resp. expanded) sub-bundle.
Then $E^c$ has dimension larger or equal to $2$ and does not admit a finer dominated splitting.
\end{conjecture}

\end{document}